\documentclass[eqno,a4paper,11pt]{article}
\title{Degeneracy in finite time of 1D quasilinear wave equations I\hspace{-.1em}I}
 \author{Yuusuke Sugiyama\footnote{e-mail:sugiyama@ma.kagu.tus.ac.jp\ \ \ telephone number:(+81)3-3260-4271}\\
Department of Mathematics,\ Tokyo University of Science\\
Kagurazaka 1-3, Shinjuku-ku, Tokyo 162-8601, Japan}
\date{}
\usepackage{amssymb,amsmath,amsthm,mathrsfs,delarray,enumerate}%mathrsfs
\usepackage{graphics}
\theoremstyle{definition} 
\newtheorem{Def}{Deffinition}[section]

\newtheorem{lemma}[Def]{Lemma}
\newtheorem{theorem}[Def]{Theorem}
\newtheorem{remark}[Def]{Remark}

\newcommand{\N}{\mathbb{N}} %set of natural number
\newcommand{\R}{\mathbb{R}} %set of real number
\makeatletter

\@addtoreset{equation}{section}

\begin{document}
%%%%%%%%%%%%%%%%  make title and   %%%%%%%%%%%%%%%%
\maketitle %タイトルをここに出力する=
%%%%%%%%%%%%%%%%  abstract   %%%%%%%%%%%%%%%%
\begin{abstract}
We consider the large time behavior of solutions to the following nonlinear wave equation:
$\partial_{t}^2 u =  c(u)^{2}\partial^2 _x u + \lambda c(u)c'(u)(\partial_x u)^2$ with the parameter $\lambda \in [0,2]$. If $c(u(0,x))$ is bounded away from a positive constant, we can construct a local solution for smooth initial data. 
However, if $c(\cdot )$ has a zero point, then $c(u(t,x))$ can be going to zero in finite time.
When $c(u(t,x))$ is going to 0 in finite time, the equation  degenerates. We give a sufficient condition that the equation with $0\leq \lambda < 2$ degenerates in finite time.
\end{abstract}
%%%%%%%%%%%%%%%%  text start  %%%%%%%%%%%%%%%%
%改行時に字下げしない。
\section{Introduction}
In this paper, we consider the Cauchy problem of the following  quasilinear wave equation: 
\begin{eqnarray} 
\left\{  \begin{array}{ll} \label{req}
  \partial_{t}^2 u =  c(u)^{2}\partial^2 _x u + \lambda c(u)c'(u)(\partial_x u)^2,  \ \ (t,x) \in (0,T] \times \R, \\   
  u(0,x) = u_0 (x),\ \ x \in \R,                \\            
    \partial_t u(0,x) = u_1 (x), \ \ x \in \R,  
\end{array} \right.  
\end{eqnarray} 
where $u(t,x)$ is an unknown real valued function, $0\leq \lambda\leq 2$ and  $c'(\theta)=dc(\theta)/d\theta$. This parameterized equation has been introduced by Glassey, Hunter and Zheng \cite{GHZ2}.
As is explained later, this equation has different  mathematical and physical backgrounds  depending on $\lambda$ (see also Chen and Shen \cite{CS}).

Throughout this paper, we assume that $c \in C^{\infty}((-1, \infty ))\cap C([-1,\infty))$ satisfies that
\begin{eqnarray} 
&&  c(\theta )>0 \ \  \mbox{for all}  \ \theta>-1, \label{c1}\\
&& c(-1)=0  , \label{c2}\\
&& c'(\theta )> 0  \ \  \mbox{for all}  \ \theta>-1 . \label{c3}
\end{eqnarray} 
Furthermore, we  assume that there exists a constant $c_0>0$ such that 
\begin{eqnarray}\label{non-deg}
c(u_0 (x))  \geq c_0
\end{eqnarray}
for all $x \in \R$.
A typical example of $c(\theta)$ satisfying \eqref{c1}-\eqref{c3} is $c(\theta)=(1+\theta)^a$ with $a>0$.

The assumptions \eqref{c1} and  \eqref{non-deg} enable us to regard the equation in \eqref{req} as a strictly hyperbolic equation near $t=0$.
By the standard local existence theorem  for strictly hyperbolic equations, the local solution of \eqref{req} with smooth initial data uniquely exists until the one of the following
two phenomena occurs. The first one is  the blow-up:
\begin{eqnarray*} 
\varlimsup_{t\nearrow T^{*}} \left( \| \partial_t u(t)\|_{L^{\infty}}+\|\partial_x u(t) \|_{L^{\infty}} \right)=\infty.
\end{eqnarray*} 
The second is the degeneracy of the equation:
\begin{eqnarray*} 
\lim_{t\nearrow T^{*}} \inf_{(s,x) \in [0,t]\times \mathbb{R}} c(u(s,x))=0.
\end{eqnarray*} 
When the equation degenerates, the standard local well-posedness theorem does not work since the equation  loses the strict hyperbolicity.
In general, for non-strictly hyperbolic equations, the persistence of the regularity of solutions does not hold (see Remark \ref{weakly}).
The aim of this paper is to give a  sufficient condition for the  occurrence of the degeneracy of the equation with $0\leq \lambda < 2$.  
The main theorem of this paper is the following.

\begin{theorem}\label{main}
Let $0\leq \lambda <2$, $(u_0, u_1)\in H^{2}(\R) \times H^1 (\R)$ and $u_1 \not\equiv 0$.
Suppose that the initial data $(u_0, u_1)$ and $c$ satisfy that \eqref{c1}-\eqref{non-deg} and
\begin{eqnarray} 
u_1 (x) \pm c(u_0 (x)) \partial_x u_0 (x) \leq 0 \ \ \ \mbox{for\ all} \ \ x \in \R. \label{inicon22} 
\end{eqnarray} 
 Then there exists  $T^{*}>0$ such that a  solution of \eqref{req}  exists uniquely and satisfies
 that $u \in C([0,T^*);H^{2}(\R))\cap C^1 ([0,T^*);H^1 (\R))$ and
\begin{eqnarray}\label{zero}
\lim_{t\nearrow T^{*}} \inf_{(s,x) \in [0,t]\times \mathbb{R}} c(u(s,x))=0.
\end{eqnarray}
Furthermore, if $0< \lambda <2$, then 
\begin{eqnarray}\label{zero2}
\lim_{t\nearrow T^{*}} c(u(t,x_0))=0
\end{eqnarray}
for some $x_0 \in \R$.
\end{theorem}

If $\lambda=2$, then the equation in \eqref{req} is formally equivalent to the following conservation system:
\begin{eqnarray*}
\partial_t \left( \begin{array}{cc} U \\ V \\ \end{array} \right)  -\partial_x \left( \begin{array}{cc} V \\ \int_{-1} ^U c(\theta)^2 d\theta \\ \end{array} \right) =0,
\end{eqnarray*}
where $U(t,x)=u(t,x)$, $V(t,x)=\int_{-\infty} ^x \partial_t u(t,y) dy$. 
This conservation system is called p-system and describes several phenomena of the wave propagation in nonlinear media  including the electromagnetic wave in a transmission line, shearing-motion in elastic-plastic rods and $1$ dimensional gas dynamics (see  Ames and Lohner \cite{al} and  Zabusky \cite{z}).
In addition to the assumptions of Theorem \ref{main}, if $\int_{\mathbb{R}} u_1 (x) dx > -2\int_{-1} ^0 c(\theta) d\theta$ is assumed, then \eqref{req} with $\lambda=2$ has 
a global smooth solution such that  the equation does not degenerate (e.g. Johnson \cite{JJ} and Yamaguchi and Nishida \cite{y-n}).
On the other hand, in \cite{s3, s4} (see Remark 1.5 in \cite{s3} and Theorem 4.1 in \cite{s4}), the author has shown  that the degeneracy \eqref{zero2} occurs in finite time,  if $\int_{\mathbb{R}} u_1 (x) dx < -2\int_{-1} ^0 c(\theta) d\theta$.
Namely, these results  say that  $-2\int_{-1} ^0 c(\theta) d\theta$ is a threshold of $ \int_{\mathbb{R}} u_1 (x) dx$ separating the global existence of solutions (such that the equation does not degenerate) and the degeneracy of the equation under the assumption \eqref{inicon22}.
If \eqref{inicon22} is not satisfied, then solutions can blow up in finite time (e.g. Klainerman and Majda \cite{km}, Manfrin \cite{mm} and Zabusky \cite{z}).
The main theorem of this paper implies that the degeneracy in finite time of the equation in \eqref{req} can occurs regardless of $\int_{\mathbb{R}} u_1 (x) dx$, when $0\leq \lambda  <2$.

When $\lambda=1$, the equation in \eqref{req} is called variational wave equation.
As its name suggests, the equation with $\lambda=1$ has a variational structure.
The variational wave equation has some physical backgrounds including nematic liquid crystal
 and long waves on a dipole chain in the continuum limit (see \cite{GHZ2}).
In \cite{GHZ, GHZ2}, Glassey, Hunter and Zheng have shown that  solutions can blow up in finite time, if \eqref{inicon22} is not satisfied (see also Remark \ref{open}).
There are a lot of papers devoted to the global existence of weak solutions to variational wave equations (e.g. Bressan and Zheng \cite{bz} and Zhang and Zheng \cite{zz1, zz2, zz3}).

When $\lambda=0$, the equation in \eqref{req} describes the wave of entropy in superfluids (e.g. Landau and Lifshitz \cite{ll}). This equation is the one dimensional version of
$$
\partial^2 _t u =c(u)^{2} \Delta u, \ \ (t,x) \in (0,T] \times \R^3,
$$
which has been studied in Lindbald \cite{lin}. 
In \cite{lin}, Lindblad has shown that solutions exists globally in time with small initial data.

In \cite{ks2, s1},  Kato and the author have shown that the equation in \eqref{req} with $c(\theta)=1+\theta$ and $\lambda=0, 1$  degenerates in finite time,
if initial data are smooth, compactly supported and satisfy \eqref{non-deg} and \eqref{inicon22}.
The main theorem of this paper removes the compactness condition on initial data and extends the result in \cite{ks2, s1} to \eqref{req} with  more general $c(\theta)$  and $0 \leq \lambda<2$.
In \cite{s1, s3}, the generalization on $\lambda$ has already been pointed out without a proof.
In fact, applying the method in \cite{s1, s3} to the equation in \eqref{req}, we can generalize the result in \cite{s1, s3} to \eqref{req} with  $0 \leq \lambda<2$ and $c(\theta)=1+\theta$.
However the compactness condition plays a crucial role in \cite{ks2, s1},
since we use the following estimates for bounded
solutions under the assumption that initial data are compactly supported:
\begin{eqnarray} \label{pre-1}
C(1+t) \geq \int_{\R} u(t,x) dx
\end{eqnarray}
and 
\begin{eqnarray}\label{pre-2}
C(1+t) \geq (2-\lambda)\int_0 ^t \int_0 ^s \int_{\R} c(u)c'(u)(\partial_x u(\tau,x))^2 dx d\tau ds.
\end{eqnarray}
\eqref{pre-1} is formally shown by the finiteness of the propagation speed and the boundedness of solutions. \eqref{pre-2} is shown by
the equation in \eqref{req}, \eqref{pre-1} and the integration by parts.
Furthermore, taking  $c(u)=1+u$, we use the following estimate in \cite{s1, s3}: 
\begin{eqnarray} \label{pre-3}
-\int_{\R} u (t,x) dx \leq C(1+t) \left( (1+t)\int_{\R} (1+u)(\partial_x u(\tau,x))^2 dx \right)^{\frac{1}{2}},
\end{eqnarray}
which is shown by  the fundamental theorem of calculus and the finiteness of the propagation speed.
We can not obtain the above estimates directly, if initial data are not compactly supported.
The first idea of the proof of Theorem \ref{main} is the use of Riemann invariant, which is a major tool for the study of $2 \times 2$ conservation systems.
When we use  the Riemann invariant in the reduction from \eqref{req} to a first order system, $\int_{-\infty} ^x c(u)c'(u)(\partial_x u(t,x))^2 dx$ appears as a force term in the first order system (see \eqref{ww} and \eqref{vv}).
The second idea for the proof is to divide the situation into the two cases that $\int_0 ^t\int_{\R} c(u)c'(u)(\partial_x u(s,x))^2 dxds$ is bounded or not.
If $\int_0 ^t\int_{\R} c(u)c'(u)(\partial_x u(s,x))^2 dxds$ is bounded, then \eqref{pre-2} holds and \eqref{pre-1} can be shown by the use of the Riemann invariant.
Hence we can use a variation of the method in \cite{ks2, s1}. 
The key for the generalization on $c(\cdot )$ is the use of $\tilde{G}(u)=\int_{-1} ^u \sqrt{c(\theta)c'(\theta)} d\theta$.
We use $\tilde{G}(u)$ in order to generalize the estimate \eqref{pre-3}.
If $\int_0 ^t\int_{\R} c(u)c'(u)(\partial_x u(s,x))^2 dxds$ is not bounded, we can use the method in \cite{s3}.
In the case that $\int_{\R} u_1(x )dx  \not\in L^1 (\R)$, the Riemann invariant can not be defined in general. Theorem \ref{main} with $u_1 \not\in L^1 (\R)$
can be shown by applying the same argument as in \cite{s3}. 

\begin{remark}
Addition to the assumptions of Theorem \ref{main}, if initial data are compactly suppurated,
then \eqref{zero2} holds for \eqref{req} with $0 \leq \lambda <2$, 
which can be shown by the finiteness of the propagation speed (see \cite{ks2, s1}).
Our method does not work for the case that $\lambda=0$ (see Remark \ref{lam0}).
\end{remark}

\begin{remark}
Under the assumptions \eqref{c1}-\eqref{non-deg}, there is still no global existence result of \eqref{req} for $0 \leq \lambda <2$. 
In stead of the assumptions, we assume that $c\in C^{\infty} (\R)$ satisfies that 
\begin{eqnarray*} 
&&  c_1 \leq c(\theta)\leq c_2 \ \ \mbox{for all}  \ \theta \in \R,\\
&& c'(\theta ) \geq0   \ \  \mbox{for all}  \ \ \theta \in \R
\end{eqnarray*} 
for some positive constants $c_1$ and $c_2$.  
Under these assumptions and  \eqref{inicon22}, Zhang and Zheng \cite{zz1} have shown that \eqref{req} has global smooth solutions with $\lambda=1$.
This global existence result has been extended to $0\leq \lambda \leq2$ in the author's paper \cite{s1}.
\end{remark}

\begin{remark}\label{open}
Here we collect some open questions for \eqref{req}.
When $\lambda=1$ and $2$, it is known that blow-up solutions exist, if \eqref{inicon22} is not satisfied.
From the first and the second equations in \eqref{fs}, we can expect that blow-up solutions exist for $0 < \lambda \leq 2$, if \eqref{inicon22} is not satisfied,
since  the right hand sides of the first and the second equations in \eqref{fs} contain $\lambda R^2$ and $\lambda S^2$ respectively, which seems to derive the singularity formation.
However,  the existence of the blow-up solution is still open, since the proofs of the blow-up theorems for $\lambda=1$ and $2$ rely on structures of the equation.
When $\lambda=0$, if $c(\cdot )$ is uniformly positive, then it seems possible  that \eqref{req} has global smooth solution for any smooth initial data, although a complete proof or a counterexample for this problem is also open.
\end{remark}

\begin{remark}\label{weakly}
It is known that a loss of the regularity appears for solutions to the following non-strictly hyperbolic equation:
$$
\partial^2 _t u - t^{2l}\partial^2 _x u - h t^{l-1}\partial_x u =0, \ \ (t,x) \in (0,\infty)\times \R,
$$
where $h$ is a constant and $l \in \N$.
Namely, in general, $(u,\partial_t u)$ does not belong to $C^1 ([0,\infty),H^s (\R)) \times C([0,\infty),H^{s-1}(\R))$ with $(u(0,x), \partial_t u(0,x)) \in H^s (\R)\times H^{s-1}(\R)$ (see Taniguchi and Tozaki \cite{TN}, Yagdjian \cite{Yad} and Qi \cite{Qi}). From this fact, we can expect that solutions of \eqref{req}  have a singularity when the equation degenerates.
\end{remark}

This paper is organized as follows: In Section 2, we recall the local well-posedness and some properties of solutions of (\ref{req}). In Sections 3 and 4, we show Theorem \ref{main} in the cases that $u_1 \in L^1 (\R)$ and $u_1 \not\in L^1 (\R) $ respectively.

\noindent
{\bf Notation}

We denote   Lebesgue space for $1\leq p\leq \infty$ and $L^2$ Sobolev space with the order $m \in \N$ on $\R$ by $L^p (\R)$  and $H^m (\R)$.
For a Banach space $X$, $C^j ([0,T];X)$  denotes the set of functions $f:[0,T] \rightarrow X$ such that $f(t)$ and its $k$ times derivatives for $k=1,2,\ldots , j$ are continuous. 
Various positive constants are simply denoted by $C$.

\section{Preliminary}
We recall the local well-posedness of \eqref{req} and some properties of solutions of \eqref{req}.
By  applying the well-known local well-posedness Theorem (e.g Hughes, Kato and Marsden \cite{HKM}, Majda \cite{m} or Taylor \cite{tay1}), we can obtain the following
theorem.
\begin{theorem} \label{lwp2} 
Let $\lambda \in \R$.
Suppose that $(u_0, u_1) \in H^{2}(\R)\times H^1 (\R)$ and that \eqref{c1} and \eqref{non-deg} hold. Then there exist $T>0$ and  a unique solution $u$ of   \eqref{req}  with
\begin{eqnarray}\label{cls21}
 u\in \bigcap_{j=0,1,2} C^{j}([0, T];H^{2-j} (\R))
\end{eqnarray}
and
\begin{eqnarray}\label{cls22}
c(u(t,x)) \geq \delta (T)\ \ \mbox{for} \ \  (t,x) \in [0,T]\times \R,
\end{eqnarray}
where $\delta (T)$ is a positive monotone decreasing function of $T \in [0,\infty)$. 
Furthermore, if $(\ref{req})$ does not have a global solution $u$ satisfying $(\ref{cls21})$ and $(\ref{cls22})$,\ then the solution $u$ satisfies 
\begin{eqnarray*} 
\varlimsup_{t\nearrow T^{*}} \left( \| \partial_t u(t)\|_{L^{\infty}}+\|\partial_x u(t) \|_{L^{\infty}} \right)=\infty
\end{eqnarray*}
 
or 
\begin{eqnarray*} 
\lim_{t\nearrow T^{*}} \inf_{(s,x) \in [0,t]\times \R} c(u(s,x))=0
\end{eqnarray*}
for some $T^{*}>0$.
\end{theorem}
We denote the maximal existence time of the solution $u$ of  (\ref{req}) constructed in Theorem \ref{lwp2} by $T^{*}$,
that is, 
\begin{align*}
T^{*} =&\sup \{ \ T> 0 \ | \ \sup_{[0,T]} \{ \| \partial_t u(t)\|_{L^{\infty}}+\|\partial_x u(t) \|_{L^{\infty}} \} <\infty,\\
&\inf_{[0,T]\times \R} c(u(t,x)) >0 \ \}.
\end{align*}

We set $R(t,x)$ and $S(t,x)$ as follows
\begin{eqnarray}\label{ri}\left\{
\begin{array}{ll}
R =\partial_t u +c(u) \partial_x u, \\
S =\partial_t u -c(u) \partial_x u.
\end{array}\right.
\end{eqnarray}
The functions $R$ and $S$ have been used in Glassey, Hunter and Zheng \cite{GHZ, GHZ2} and Zhang and Zheng \cite{zz1}.
We recall some properties of $R$ and $S$ proved in \cite{s1}.  

By (\ref{req}), $R$ and $S$ are  solutions to the system of the following first order equations:
\begin{eqnarray}\label{fs}
 \left\{  
\begin{array}{ll}
\partial_t R -c(u)\partial_x R=\dfrac{c'(u)}{2c(u)}(R S-S ^2) 
+\lambda  \dfrac{c'(u)}{4c(u)}(R -S )^2 , \\
\partial_x u = \dfrac{1}{2c(u)}(R - S),\\
\partial_t S +c(u)\partial_x S =\dfrac{c'(u)}{2c(u)}(S R-R ^2)
+\lambda  \dfrac{c'(u)}{4c(u)}(S -R )^2 . 
\end{array}
\right.
\end{eqnarray} 

\begin{lemma}\label{zz1} 
Let $0 \leq \lambda \leq 2$. Suppose that the assumptions of Theorem \eqref{req} are satisfied. Then  we  have
\begin{eqnarray}\label{lemmazzes1}
R(t,x),\ S(t,x) \leq 0 \ \mbox{for} \ (t,x) \in [0,T ^* ) \times \R,
\end{eqnarray}
where $R$ and $S$ are the functions in $(\ref{ri})$ for the solution $u$   of \eqref{req} constructed by Theorem \ref{lwp2}.
\end{lemma}

\begin{lemma}\label{zz2} 
Let $p\geq \max\{2,\frac{2}{\lambda}\}$. Suppose that the assumptions of Theorem \eqref{req} are satisfied. Then  we  have for $0 < \lambda \leq 2$
\begin{eqnarray}\label{lemmazzes2} 
\| R (t) \|^p _{L^{p}}+\| S (t)\|^p _{L^{p}}\leq \|R (0) \|^p _{L^{p}}+\| S (0) \|^p _{L^{p}}, \ \mbox{for} \ t \in [0,T^* ) ,
\end{eqnarray}
where $R$ and $S$ are the functions in $(\ref{ri})$ for the solution $u$  of \eqref{req} constructed by Theorem \ref{lwp2}.
Furthermore $\displaystyle \|R(t)\|_{L^{\infty}}$ and $\displaystyle \|S(t)\|_{L^{\infty}}$ are uniformly bounded with $t \in [0.T^* )$ for $0 \leq \lambda \leq 2$. 
\end{lemma}

Lemmas \ref{zz1} and \ref{zz2} have been shown in the author's paper \cite{s1}.
The proofs are essentially the same  as in the case that $\lambda  =1$, which are proved in  Zhang and  Zheng \cite{zz1}.
In \cite{s1, s3}, it is assumed only $p\geq 2/\lambda$ for the inequality \eqref{lemmazzes2}. But the proof in \cite{s1} is not collect  for 
$p<2$. In fact, in \cite{s1}, the proof of \eqref{lemmazzes2} is based on the following inequality:
$$
\frac{1}{p}\frac{d}{dt}  \int_{\R}  \tilde{R}^p  +  \tilde{S}^p dx \leq  -(\frac{1}{2} -\frac{\lambda }{4})\int_{\R} \frac{c'(u)}{c(u)}  \Tilde{R} \Tilde{S} ( \Tilde{R}-\Tilde{S})(( \Tilde{R})^{p-2} -(\Tilde{S})^{p-2}) dx,
$$
where $\Tilde{R} = -R$ and $\Tilde{S}=-S$.
If $p<2$, the right hand side of this inequality is not  negative except for $\lambda=2$.
However, in \cite{s1}, we only use \eqref{lemmazzes2} for $p \geq 2$.
If $\lambda=2$, \eqref{lemmazzes2} holds for all $p\geq 1$.

We note that \eqref{lemmazzes1} implies that $\partial_t u(t,x) \leq 0$ for all $(t,x) \in [0,T ^* ) \times \R$.

\section{Proof of Theorem \ref{main} with $u_1 \in L^1 (\R)$}

We show Theorem \ref{main} in the case that  $u_1 \in L^1 (\R)$.

First, we show  \eqref{zero} by the contradiction argument. From Theorem \ref{lwp2} and Lemma \ref{zz2}, it is enough to show that $T^* < \infty$.
We set $G(u)=\int_{-1} ^u c(\theta) d\theta$ for $u \geq -1$ and $\mu=2-\lambda$ and define the Riemann invariants $(w_1(t,x) , w_2 (t,x))$ and $(v_1(t,x) , v_2 (t,x))$ as follows:
\begin{eqnarray*}
\begin{array}{ll}
w_1 = \displaystyle\int_{-\infty} ^x \partial_t u dx +G(u), \\
w_2 = \displaystyle\int_{-\infty} ^x \partial_t u dx -G(u)
\end{array}
\end{eqnarray*}
and
\begin{eqnarray*}
\begin{array}{ll}
v_1 = \displaystyle\int_{x} ^\infty \partial_t u dx -G(u), \\
v_2  = \displaystyle\int_{x} ^\infty \partial_t u dx +G(u).
\end{array}
\end{eqnarray*}
From \eqref{req}, $(w_1(t,x) , w_2 (t,x))$ and $(v_1(t,x) , v_2 (t,x))$ satisfy that
the following systems:
\begin{eqnarray}\label{ww}\left\{
\begin{array}{ll} 
\partial_t w_1 -c(u) \partial_x w_1 =-\mu \int_{-\infty} ^x  \tilde{e} (t,y) dy,\\
  \partial_t w_2 +c(u) \partial_x w_2 =-\mu \int_{-\infty} ^x   \tilde{e} (t,y) dy
\end{array}\right.
\end{eqnarray}
and
\begin{eqnarray}\label{vv}\left\{
\begin{array}{ll}
 \partial_t v_1 -c(u) \partial_x v_1 =-\mu \int_{x} ^{\infty}  \tilde{e} (t,y) dy,\\
 \partial_t v_2 +c(u) \partial_x v_2 =-\mu \int_{x} ^{\infty}   \tilde{e} (t,y) dy,
\end{array}\right.
\end{eqnarray}
where $ \tilde{e} (t,y) = c'(u)c(u) (\partial_x u )^2 (t,y)$.
Let $x_{\pm} (t)$ be characteristic curves on the first and third equations of \eqref{ww}
respectively. That is,  $x_{\pm} (t)$ are solutions to the following differential equations:
 \begin{eqnarray}\label{cc}
\dfrac{d}{dt} x_{\pm} (t)=\pm c(u(t,x_{\pm} (t))).
\end{eqnarray}
\eqref{ww} and \eqref{vv} imply that
\begin{align}
w_1 (t,x_{-} (t))= w_1 (t_0, x_- (t_0)) -\mu \int_{t_0} ^t \int_{-\infty} ^{x_{-} (s)} \tilde{e}(s,y) dy ds \label{ww10} 
\end{align}
and
\begin{align}
w_2 (t,x_{+} (t))= w_2 (t_0, x_+ (t_0)) -\mu \int_{t_0} ^t \int_{-\infty} ^{x_{+} (s)} \tilde{e}(s,y) dy ds  \label{ww20}.
\end{align}

{\bf Case that $\int_0 ^t \int_{\R}   \tilde{e}(s,y) dy ds$ is bounded.}

By the contradiction argument, we show that $T^*$ is finite in the case that $\int_0 ^t \int_{\R}   \tilde{e}(s,y) dy ds$ is bounded on $[0,\infty)$.
We suppose that $T^* = \infty$.
\eqref{ww10} and \eqref{ww20} imply that
\begin{align}
w_1 (t,x_{-} (t))\geq  w_1 (0, x_{-}(0)) -\mu \int_{0} ^{\infty} \int_{-\infty} ^{x_{-} (s)} \tilde{e}(s,y) dy ds \label{ww1} 
\end{align}
and
\begin{align}
w_2 (t,x_{+} (t))\geq  w_2 (0, x_{+}(0)) -\mu \int_{0} ^{\infty} \int_{-\infty} ^{x_{+} (s)} \tilde{e}(s,y) dy ds. \label{ww2} 
\end{align}

We fix an arbitrary number  $\varepsilon >0$.
Since $ \lim_{|x| \rightarrow  \infty}u_0 (x) =0$, $u_1 \in L^1 (\R)$ and $\partial_x w_j (0,x) \leq 0$ for $j=1, 2$, there exists a constant $M_0>0$ such that
$$
G(0)- \varepsilon \leq w_1 (0,x)\leq G(0)
$$
and
$$
-G(0)-\varepsilon \leq w_2 (0,x)\leq -G(0)
$$
for any $x \leq -M_0$. Noting $  x_{\pm} (t)$ goes to $-\infty$  as $x_{\pm} (0) \rightarrow -\infty $ for all $t \geq 0$,
since $\int_0 ^\infty \int_{\R}   \tilde{e}(s,y) dy ds$ is bounded, we have by the Lebesgue convergence theorem
$$
\lim_{x_{\pm} (0) \rightarrow -\infty} \int_{0} ^{\infty} \int_{-\infty} ^{x_{\pm} (s)} \tilde{e}(s,y) dy ds =0.
$$
Hence, from \eqref{ww1} and \eqref{ww2}, there exists  a constant $M_1>0$  such that
if $x_{\pm} (0) \leq -\max\{M_0, M_1\}$, then for all $t \geq 0$
\begin{eqnarray*}
G(0)-\varepsilon \leq w_1 (t,x_{-} (t)) \leq  G(0) 
\end{eqnarray*}
and 
\begin{eqnarray*}
-G(0)-\varepsilon \leq w_2 (t,x_{+} (t)) \leq  -G(0).
\end{eqnarray*}
We note the positive constant $M_1$ can be chosen independently of $t$.
Hence the equality $2G(u(t,x))=w_1 (t,x)-w_2 (t,x)$ yields that
$$
G(0)- \varepsilon \leq G(u(t, x)) \leq G(0) +\varepsilon,
$$
if $x\leq x_{-}(t)$, where $x_{-}(0)\leq -\max\{M_0, M_1\}$.
Since $G$ is invertible and $G^{-1}$ is continuous, this inequality implies that
\begin{eqnarray} \label{u-ip}
|u(t,x)|\leq C\varepsilon
\end{eqnarray}
 with $x\leq x_{-}(t)$, where $x_{-}(0)\leq -\max\{M_0, M_1\}$.
From the above estimates of $w_1$ and $G(u)$, we have
$$
   -C\varepsilon \leq  \int_{-\infty} ^{x_{-} (t)} \partial_t u(t,y) dy \leq 0.
$$
Since 
\begin{align*}
\dfrac{d}{dt} \int_{-\infty} ^{x_{-} (t)} u(t,x) -u_0 (x)dx = & \int_{-\infty} ^{x_{-} (t)} \partial_t u(t,x) dx \\
& - (u(t,x_{-} (t))-u_0 (x_{-} (t)))c(u(t,x_{-}(t))),
\end{align*}
 if $x_{-}(0)\leq -\max\{M_0, M_1\}$, then we have from \eqref{u-ip}
\begin{eqnarray}\label{es-mi0}
 -\int_{-\infty} ^{x_{-} (t)} u(t,y)-u_0 (x) dy\leq C\varepsilon t .
\end{eqnarray}
By using \eqref{vv}, we have in the same way as in the derivation of \eqref{es-mi0}
\begin{eqnarray}\label{es-mi2}
-\int_{x_+ (t)} ^{\infty} u(t,y)-u_0 (x) dy\leq C\varepsilon t ,
\end{eqnarray}
if $x_{+}(0) \geq  M_2$ for  sufficiently large $M_2 >0$.

We set $F(t)=-\int_{\R} u(t,x)- u_0 (x)dx$ and take $M \geq \max\{M_0, M_1, M_2 \}$.
From the integration by parts and $\eqref{req}$, it follows that
$F''(t)=\mu \int_{\R} \tilde{e}(t,x) dx\geq 0$. Integrating this equality twice on $[0,t]$ and dividing by $t$,  we have
\begin{eqnarray*}
\frac{F(t)}{t} \geq F'(0).
\end{eqnarray*}
By \eqref{es-mi} and \eqref{es-mi2}, we have 
\begin{align}
F'(0) \leq \frac{F(t)}{t} = &  \frac{-1}{t}\left(\int_{-\infty} ^{x_- (t)}+\int_{x_{-}(t)} ^{x_+ (t)}  +\int_{x_+ (t)} ^\infty  \right)u(t,x)-u_0 (x) dx \notag \\
\leq &    C\varepsilon  - \frac{1}{t} \int_{x_{-}(t)} ^{x_+ (t)} u(t,x)-u_0 (x) dx. \label{fff2}
\end{align}
Now we estimate the second term of the right hand side of \eqref{fff2}. We set $\tilde{G}(u)=\int_{-1} ^u \sqrt{c(\theta)c'(\theta)} d\theta$.
The Schwarz inequality implies that
$$
\tilde{G}(u)^2\leq \int_{-1} ^u c'(\theta) d\theta \int_{-1} ^{u} c(\theta) d\theta =c(u)\int_{-1} ^{u} c(\theta) d\theta.
$$
Hence $\tilde{G}(u)$ can be defined for $u \geq -1$. From \eqref{c1} and \eqref{c2}, we have that $\tilde{G}'(u)= \sqrt{c(u)c'(u)} >0$, from which $\tilde{G}(\cdot )$
is invertible on $[0,\infty)$ and $\tilde{G}^{-1}(\cdot )$ is continuous.
The fundamental theorem of calculus yields that
\begin{align} \label{qq}
\tilde{G} (u(t,x))=\tilde{G} (u(t,x_{-}(t))) + \int_{x_{-} (t)} ^x \sqrt{c(u)c'(u)} \partial_y u(t,y) dy.
\end{align}
Applying \eqref{u-ip} to the first term of the right hand side of \eqref{qq} and the Schwarz inequality to the second term, we have
\begin{eqnarray*}
\tilde{G} (u(t,x)) \geq \tilde{G} (0)-C\varepsilon - \sqrt{|x_{+}(t)-x_{-}(t)| \int_{\R} \tilde{e}(t,y) dy}
\end{eqnarray*}
for $x \in [x_{-}(t), x_+ (t)]$. 
Since $\tilde{G}^{-1}$ is a monotone increasing function,  we have 
\begin{eqnarray*}
u(t,x) \geq \tilde{G}^{-1} \left( \tilde{G} (0)-C\varepsilon - \sqrt{|x_{+}(t)-x_{-}(t)| \int_{\R} \tilde{e}(t,y) dy}\right).
\end{eqnarray*}
From \eqref{cc} and \eqref{u-ip}, we have $|x_{+}(t)-x_{-}(t)| \leq C_M + C^* t$, where $C_M >0$ depends on $M_j$ for $j=1,2,3$ and $C^* >0$ can be chosen independently of the  three constants. 
Integrating the both sides of this inequality on $[x_{-}(t), x_+ (t)]$,
we have
\begin{align} 
-\int_{x_{-}(t)} ^{x_+ (t)} u(t,x) dx \leq &-(C_M + C^* t) \tilde{G}^{-1} \Biggl( \tilde{G} (0)-C\varepsilon  \notag \\
&  - \sqrt{(C_M + C^* t) \int_{\R} \tilde{e}(t,y) dy}\Biggr). \label{u-int}
\end{align}
While the Schwarz inequality implies that
\begin{eqnarray*}
\int_{x_{-}(t)} ^{x_+ (t)} |u_0 (x)| dx \leq C\sqrt{C_M + C^* t}\|u_0 \|_{L^2}.
\end{eqnarray*}
From this inequality, \eqref{fff2} and \eqref{u-int}, we have
\begin{align}
F'(0) \leq & \frac{C\sqrt{C_M + C^* t}\|u_0 \|_{L^2}}{t} + C\varepsilon \notag \\
&- \frac{C_M + C^* t}{t} \tilde{G}^{-1} \left( \tilde{G} (0)-C\varepsilon - \sqrt{(C_M + C^* t) \int_{\R} \tilde{e}(t,y) dy}\right).\label{las-case1}
\end{align}
Since we assume that $ \int_0 ^{\infty} \int_{\R} \tilde{e}(s,y)dyds < \infty$, there exists a monotone increasing sequence $\{t_j \}_{j \in \N}$
such that $ \lim_{j \rightarrow \infty} t_j = \infty$ and
$$
\lim_{j\rightarrow \infty} (C_M + C^* t_j )\int_{\R}\tilde{e}(t_j ,y)dy =0.
$$
Putting $t=t_j$  in \eqref{las-case1} and  taking $j \rightarrow \infty$, since $\tilde{G}$ is continuous, we obtain
\begin{eqnarray*}
F'(0) \leq  C \varepsilon -C^*\tilde{G}^{-1} \left( \tilde{G} (0)-C\varepsilon \right)  \leq C \varepsilon,
\end{eqnarray*}
which contradicts to the assumption that $u_1 \not\equiv 0$, if $\varepsilon $ is sufficiently small.
Therefore we obtain $T^{*} < \infty$ in the case that $\int_0 ^t \int_{\R}   \tilde{e}(s,y) dy ds$ is bounded.

\begin{remark}
The strictly positivity of $c'$ is only used in the estimate of $\tilde{G}$.
It is enough to assume that $c'(\theta) \geq 0$ for $\theta >0$ in the case that $\int_0 ^t \int_{\R}   \tilde{e}(s,y) dy ds$ is unbounded.
In this case, we use the method in \cite{s3, s4}.
In \cite{s3, s4}, instead of \eqref{c3}, it is assumed that $c'(\theta)\geq 0$ for the occurrence of the degeneracy of the equation in \eqref{req} with $\lambda=2$.
\end{remark}

{\bf Case that $\int_0 ^t \int_{\R}   \tilde{e}(s,y) dy ds$ is unbounded.}

We suppose that $T^* = \infty$.
In this case, from the identity $$ -\int_{\R} \partial_t u(t,x) dx =-\int_{\R} u_1 (x) dx +\mu \int_0 ^t \int_{\R} \tilde{e}(s,x)dx ds,$$
there exists a positive number $T$ such that 
\begin{eqnarray}\label{kar-asmp}
-\int_{\R} \partial_t u(T,x) dx > 2G(0).
\end{eqnarray}
From \eqref{ww1} and \eqref{ww2}, if the plus and minus characteristic curves cross at some point $(t_0, x_0)$ with $t_0 \geq T$, then   
we have  by the definitions of $w_1$ and $w_2$
\begin{align}\label{es-u}
2G(u(t_0,x_0))=& w_1 (T,x_{-} (T)) -w_2 (T,x_{+} (T))- \mu \int_{T} ^{t_0} \int_{x_{+} (s)} ^{x_{-} (s)} \tilde{e}(s,y) dy ds\notag \\
 =& \displaystyle\int_{x_{+} (T)} ^{x_{-}(T)} \partial_t u (T, x) dx+ G(u (T,x_{+} (T))) +G(u (T, x_{-} (T)))\notag \\
 & -\mu \int_{T} ^{t_0} \int_{x_{+} (s)} ^{x_{-} (s)} \tilde{e}(s,y) dy ds.
\end{align}
By \eqref{kar-asmp} and the facts that $\lim_{|x| \rightarrow \infty} u(T,x)=0$ and that $\partial_t u(T,\cdot ) \in L^1 (\R)$, there exists a number $M>0$ such that
\begin{eqnarray}\label{kar-asmp2}
\int_{-M} ^{M} \partial_t u (T, x) dx+ G(u(T,-M))+G(u(T, M)))<0.
\end{eqnarray} 
We set $\Tilde{F}(t)=-\int_{\R} u(t,x)- u (T, x)dx$.
We derive a estimate of $\Tilde{F}(t)$ which contradicts to \eqref{kar-asmp2}.

Suppose that the plus and minus characteristic curves $x_{\pm} (t)$ defined in \eqref{cc} pass through $(T, \mp M)$ respectively.
The characteristics $x_{\pm}(t)$ are drawn on the $(x,t)$ plane as follows:
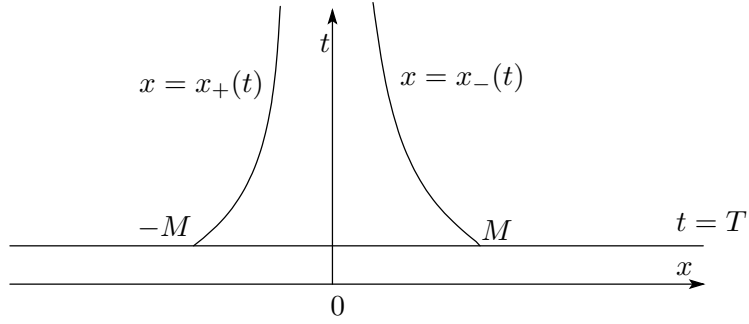
\begin{center}{Figure 1: the two characteristic curves on  the $(x,t)$ plane}
{\unitlength 0.1in%
\begin{picture}(40.6000,18.1000)(4.1000,-18.4000)%
% VECTOR 2 0 3 0 Black White  
% 2 630 1600 4220 1600
% 
\special{pn 8}%
\special{pa 630 1600}%
\special{pa 4220 1600}%
\special{fp}%
\special{sh 1}%
\special{pa 4220 1600}%
\special{pa 4153 1580}%
\special{pa 4167 1600}%
\special{pa 4153 1620}%
\special{pa 4220 1600}%
\special{fp}%

\special{pn 8}%
\special{pa 630 1400}%
\special{pa 4220 1400}%
\special{fp}%

% VECTOR 2 0 3 0 Black White  
% 2 2300 1600 2300 170
% 
\special{pn 8}%
\special{pa 2300 1600}%
\special{pa 2300 170}%
\special{fp}%
\special{sh 1}%
\special{pa 2300 170}%
\special{pa 2280 237}%
\special{pa 2300 223}%
\special{pa 2320 237}%
\special{pa 2300 170}%
\special{fp}%
% STR 2 0 3 0 Black White  
% 4 4080 1540 4080 1590 2 0 0 0
% x
\put(40.8000,-13.2000){\makebox(0,0)[lb]{$t=T$}}%
\put(40.8000,-15.5000){\makebox(0,0)[lb]{$x$}}%
% STR 2 0 3 0 Black White  
% 4 2235 300 2235 350 2 0 0 0
% t
\put(22.3500,-3.9000){\makebox(0,0)[lb]{$t$}}%
% STR 2 0 3 0 Black White  
% 4 570 530 570 580 2 0 0 0
% a
% b
\put(13.0000,-6.3500){\makebox(0,0)[lb]{$x=x_+ (t)$}}%
% STR 2 0 3 0 Black White  
% 4 2690 550 2690 600 2 0 0 0
% c
\put(26.5000,-6.0000){\makebox(0,0)[lb]{$x=x_- (t)$}}%
% STR 2 0 3 0 Black White  
% 4 3990 760 3990 810 2 0 0 0
% d
% M
\put(30.7000,-13.7000){\makebox(0,0)[lb]{$M$}}%
% STR 2 0 3 0 Black White  
% 4 1190 1710 1190 1760 2 0 0 0
% -M
\put(12.9000,-13.6000){\makebox(0,0)[lb]{$-M$}}%
% SPLINE 2 0 3 0 Black White  
% 5 3350 1600 3860 1370 4320 710 4470 440 4470 440
% 
\put(22.9000,-17.6000){\makebox(0,0)[lb]{$0$}}%
\special{pn 8}%
\special{pa 3064 1400}%
\special{pa 3044 1378}%
\special{pa 3019 1358}%
\special{pa 2995 1338}%
\special{pa 2971 1317}%
\special{pa 2948 1296}%
\special{pa 2925 1274}%
\special{pa 2881 1230}%
\special{pa 2841 1184}%
\special{pa 2822 1160}%
\special{pa 2786 1110}%
\special{pa 2754 1058}%
\special{pa 2724 1004}%
\special{pa 2711 977}%
\special{pa 2698 948}%
\special{pa 2685 920}%
\special{pa 2673 891}%
\special{pa 2662 862}%
\special{pa 2651 832}%
\special{pa 2641 802}%
\special{pa 2632 771}%
\special{pa 2623 741}%
\special{pa 2614 710}%
\special{pa 2598 646}%
\special{pa 2590 615}%
\special{pa 2583 582}%
\special{pa 2576 550}%
\special{pa 2570 517}%
\special{pa 2564 485}%
\special{pa 2558 452}%
\special{pa 2552 418}%
\special{pa 2542 352}%
\special{pa 2537 318}%
\special{pa 2532 285}%
\special{pa 2517 183}%
\special{pa 2513 149}%
\special{pa 2510 130}%
\special{fp}%
% SPLINE 2 0 3 0 Black White  
% 4 1280 1600 1760 1230 2030 150 2030 150
% 
\special{pn 8}%
\special{pa 1581 1400}%
\special{pa 1631 1360}%
\special{pa 1700 1297}%
\special{pa 1721 1275}%
\special{pa 1742 1252}%
\special{pa 1761 1229}%
\special{pa 1780 1205}%
\special{pa 1814 1155}%
\special{pa 1830 1129}%
\special{pa 1845 1103}%
\special{pa 1859 1076}%
\special{pa 1872 1048}%
\special{pa 1885 1021}%
\special{pa 1897 992}%
\special{pa 1908 963}%
\special{pa 1918 934}%
\special{pa 1928 904}%
\special{pa 1937 874}%
\special{pa 1953 812}%
\special{pa 1960 781}%
\special{pa 1967 749}%
\special{pa 1979 685}%
\special{pa 1985 652}%
\special{pa 1989 619}%
\special{pa 1994 586}%
\special{pa 1998 553}%
\special{pa 2006 485}%
\special{pa 2015 383}%
\special{pa 2017 349}%
\special{pa 2020 314}%
\special{pa 2022 279}%
\special{pa 2024 245}%
\special{pa 2028 175}%
\special{pa 2030 150}%
\special{fp}%
\end{picture}}
\end{center}
From \eqref{es-u} and \eqref{kar-asmp2}, these characteristic curves $x_{+}(t)$ and $x_{-}(t)$ do not cross for all $t\geq T$.
Hence it follows that
\begin{eqnarray}\label{lim-u}
\lim_{t\rightarrow \infty} c(u(t,x_{\pm} (t))) =0.
\end{eqnarray} 
And $\Tilde{F}$ can be divided as follows:
\begin{eqnarray} \label{F-div}
\Tilde{F}(t)=-\left(\int_{-\infty} ^{x_+ (t)}+\int_{x_{+}(t)} ^{x_- (t)}  +\int_{x_- (t)} ^\infty  \right)u(t,x)-u (T, x) dx.
\end{eqnarray}
Now we estimate $\int_{-\infty} ^{x_{+} (t)} u(t,x)- u (T,x) dx$.
From \eqref{cc}, we have 
\begin{align}
\dfrac{d}{dt} \int_{-\infty} ^{x_{+} (t)} u(t,x)-u (T,x) dx=& \int_{-\infty} ^{x_{+} (t)} \partial_t u(t,x) dx \notag \\
& + (u(t,x_{+} (t))-u (T, x_{+} (t)))c(u(t,x_{+}(t))).  \label{tsui} 
\end{align}
From \eqref{ww20} and the definition of $w_2$, the first term of the right hand side of \eqref{tsui} can be estimated as follows:
\begin{align*}
\int_{-\infty} ^{x_{+} (t)} \partial_t u(t,x) dx = & G(u(t, x_{+}(t)))+ w_2 (T,-M)-\mu  \int_T ^t \int_{-\infty} ^{x_+ (s)} \tilde{e}(s,x)dx ds\\
\geq &  w_2 (T,-M)-\mu  \int_T ^t \int_{-\infty} ^{x_+ (s)} \tilde{e}(s,x)dx ds.
\end{align*}
From the boundedness of $u$, the second term of the right hand side of \eqref{tsui} can be estimated as
$$
(u(t,x_{+} (t))-u (T, x_{+} (t)))c(u(t,x_{+}(t)))\geq -Cc(u(t,x_{+}(t))).
$$
Hence we have from \eqref{tsui}
\begin{align}
\int_{-\infty} ^{x_{+} (t)} u(t,x)-u (T,x) dx\geq & (t-T)w_2 (T,-M) -C\int_T ^t c(u(s,x_{+}(s))) ds \notag \\
&- \mu \int_T ^t \int_T ^s \int_{-\infty} ^{x_+ (\tau)} \tilde{e}(\tau,x)dx d\tau ds. \label{esf1}
\end{align}
Using $v_1$ instead of $w_2$, in the same way as in the derivation of $\eqref{esf1}$, we get
\begin{align}
\int_{x_{-} (t)} ^{\infty} u(t,x)-u (T,x) dx\geq & (t-T)v_1 (T,-M) -C\int_T ^t c(u(s,x_{-}(s)))ds \notag \\
&- \mu \int_T ^t \int_T ^s \int_{x_- (\tau)} ^{\infty} \tilde{e}(\tau,x)dx d\tau ds. \label{esf2}
\end{align}
The boundedness of $u(t,x)$ and  $|x_{+}(t)-x_- (t) |$ yields that
\begin{eqnarray} \label{es-mi}
-\int_{x_{+} (t)} ^{x_- (t)} u(t,x)-u (T,x) dx \leq C|x_- (t) - x_+ (t)| \leq C.
\end{eqnarray}
While, in the same way as in the computation of $F$,  we have
\begin{eqnarray}\label{ff}
\Tilde{F}(t)=(t-T)\Tilde{F}'(T) +\mu \int_T ^t \int_T ^s \int_\R \tilde{e}(\tau,x)dx d\tau ds.
\end{eqnarray}
By \eqref{F-div}, \eqref{esf1}, \eqref{esf2} and \eqref{es-mi}, we have 
\begin{align}
\Tilde{F}(t) 
\leq &  C-(t-T)(w_2 (T,-M)+v_1 (T,M))\notag \\
&  + C\int_T ^t   c(u(s,x_{+}(s))) +c(u(s,x_{-}(s))) ds  \notag \\
&+ \mu \int_T ^t \int_T ^s \left( \int_{-\infty} ^{x_+ (\tau)}  + \int_{x_+ (\tau)} ^{\infty} \right) \tilde{e}(\tau,x)dx d\tau ds. \label{fff}
\end{align}
From the definitions of $w_2$ and $v_1$, the second term of the right hand side of \eqref{fff} can be written as follows:
\begin{align*}
w_2 (T,-M)+ v_1 (T,M)= & \left(\int_{-\infty} ^{-M} +\int_{M} ^{\infty} \partial_t u(T,x) dx \right) \\
&-(G(u(T,M)) +G(u(T,-M))), 
\end{align*}
from which, \eqref{ff} and \eqref{fff} yield that
\begin{align*}
  - \int_{-M} ^{M} \partial_t u(T,x) dx \leq & \frac{-\mu}{(t-T)}\int_T ^t \int_T ^s \int_{x_{+}(\tau)} ^{x_- (\tau)} \tilde{e}(\tau,x)dx d\tau ds \\
 &G(u(T,M)) +G(u(T,-M))\\
&+\frac{C}{t-T}\int_T ^t   c(u(s,x_{+}(s))) +c(u(s,x_{-}(s)))ds\\
 \leq & G(u(T,M)) +G(u(T,-M)) \\
&+\frac{C}{t-T}\int_T ^t   c(u(s,x_{+}(s))) +c(u(s,x_{-}(s)))ds.
\end{align*}
From \eqref{lim-u}, the second term of the right hand side of the above inequality tends to $0$ as $t\rightarrow \infty$.
Hence, taking $t\rightarrow \infty$ in  the above inequality, we have
\begin{eqnarray*}
 - \int_{-M} ^M \partial_t u(T,x) dx\leq G(u(T,M)) +G(u(T,-M)),
\end{eqnarray*}
which contradicts to \eqref{kar-asmp2}. 
Hence we have that $T^*<\infty$  in the case  $\int_0 ^t \int_{\R}   \tilde{e}(s,y) dy ds$ is unbounded.

From the above argument of the two cases, we have $T^* < \infty$ and $\eqref{zero}$.

Next we give an outline of the proof of \eqref{zero2} for $0< \lambda <2$. 
The proof is the same as in \cite{ks2, s1, s3}.
Since $u(t,x)$ is a monotone decreasing function with $t$ for all $x \in \R$, we can define
$\tilde{u}(x)=\lim_{t\nearrow T^*} u(t,x)$. 
While, \eqref{lemmazzes2} in Lemma \ref{zz2} implies that $\| c(u) \partial_x u(t)\|_{L^p}$
is uniformly bounded with $t \in [0,T^*)$ and $p=\max\{2,2/\lambda\}$. Hence, by the standard argument on the Sobolev space,
it follows that $G(\tilde{u})-G(0) , c(\tilde{u}) \partial_x \tilde{u} \in L^p (\R)$ and that $G(\tilde{u}(\cdot ))$ is a continuous function. Therefore we have
$
\lim_{|x|\nearrow \infty } G(\tilde{u})-G(0)=0,
$
from which, the continuity of $G(\tilde{u}(\cdot ))$ implies that $ G(\tilde{u}(x_0))=0$ for some $x_0 \in \R$.
While, from the monotonicity of $G(u(t,x))$ with $t$, we have
$$
\lim_{t\nearrow T^{*}} \inf_{(s,x) \in [0,t]\times \mathbb{R}} G(u(s,x))=\inf_{x \in \R} \lim_{t\nearrow T^*} G(u(t,x))=\inf_{x \in \R} G(\tilde{u}(x)).
$$
Hence, by the continuity of $G^{-1}$, we have $\lim_{t \nearrow T^* }u(t,x_0)=-1$, which implies \eqref{zero2}.

\begin{remark}\label{lam0}
In the case that $\lambda=0$, since the boundedness of $\| c(\Tilde{u}) \partial_x \Tilde{u}(t)\|_{L^p}$ is unknown for $p \not=\infty$, the above argument does not work.
Hence the case that $\lambda=0$ is excluded in \eqref{zero2}. 
\end{remark}

\section{Proof of Theorem \ref{main} with $u_1 \not\in L^1 (\R)$}

By using the same argument as in \cite{s3}, we can show Theorem \ref{main} with $u_1 \not\in L^1 (\R)$.
We show that the degeneracy \eqref{zero} occurs in finite time for the reader's convenience.
\eqref{zero2} can be shown by same way as in the case that $u_1 \in L^1 (\R)$.
In the same argument as in Section $3$, we can say that then \eqref{zero} occurs at $T^*$, if $T^* < \infty$.
Hence it is enough to show that $T^*$ is finite.
We define a cut-off function $\psi  \in C^{\infty}_0 (\mathbb{R})$ as
\begin{eqnarray*} 
\psi (x) =\left\{ \begin{array}{ll}
 1\ \ \mbox{for} \ |x|\leq 1,\\
0\ \ \mbox{for} \ |x|\geq 2
\end{array}
\right. 
\end{eqnarray*}
and $0 \leq \psi (x) \leq 1$.
We set $ \psi_\varepsilon (x)=  \psi (\varepsilon x) $ and $ F_\varepsilon (t) = -\int_{\mathbb{R}}  \psi_\varepsilon (x) u(t,x) dx$.
From \eqref{req} and the integration by parts, we have
\begin{align*}
F'' _\varepsilon (t) =&-\varepsilon^2 \int_{\R}  \psi'' (\varepsilon x)  G_2 (u)  dx+ \mu \int_{\R}  \psi _\varepsilon (x)  c(u) c'(u)(\partial_x u)^2  dx\\
\geq &-\varepsilon^2 \int_{\R}  \psi'' (\varepsilon x) G_2 (u)   dx,
\end{align*}
where $G_2 (u)=\int_{-1} ^u c(\theta)^2 d\theta$.

Since $-1 \leq u(t,x)\leq u(0,x) \leq C$, it follows that
$$
F'' _\varepsilon (t)\geq -C\varepsilon .
$$
Namely we have
$$
 F_\varepsilon (t)\geq F_{\varepsilon }(0) + t F' _{\varepsilon }(0) -C\varepsilon t^2 .
$$
The boundedness of $u(t,x)$ with $(t,x) \in [0,\infty) \times \mathbb{R}$ yields that
$$
 |F_\varepsilon (t)| \leq \int_{\mathbb{R}} |\psi_{\varepsilon} (x)| |u(t,x)| dx \leq  \frac{C}{\varepsilon}.
$$
Hence we have
\begin{eqnarray*}
\frac{C}{\varepsilon}+C\varepsilon t^2 \geq  t F' _{\varepsilon }(0).
\end{eqnarray*}
Putting $ \varepsilon = 1/(1+t)$, we  have
\begin{eqnarray*}
\frac{2C(t+1)}{t}\geq  -\int_{\mathbb{R}} \psi \left(\frac{x}{t+1}\right) u_1 (x) dx.
\end{eqnarray*}
Since $ u_1 \not\in L^1 (\mathbb{R})$, the right hand side is going to infinity as $ t\rightarrow \infty$, which is a contradiction.
Therefore we have $T^* < \infty$.





\end{document}